\documentclass[draft]{amsart}

\usepackage{amsmath}
\usepackage{amsthm}
\usepackage{amssymb}
\usepackage[all]{xy}

\UseComputerModernTips

\newcommand{\id}[1]{\ensuremath{\operatorname{id}_{#1}}}
\newcommand{\eins}{1 \hspace{-2.3pt} \mathrm{l}}

\newcommand{\chic}{\ensuremath{\chi_c}}
\newcommand{\chicmot}{\ensuremath{\chi_{c,\mathrm{mot}}}}
\newcommand{\chimot}{\ensuremath{\chi_{\mathrm{mot}}}}
\newcommand{\tate}{\ensuremath{\mathbb{L}}}
\newcommand{\tatemot}{\ensuremath{\tate_{\mathrm{mot}}}}

\newcommand{\Knull}{\ensuremath{\operatorname{K}_0}}
\newcommand{\KnullG}{\ensuremath{\operatorname{K}_0^G}}
\newcommand{\KnullGstr}{\ensuremath{\operatorname{K^{\prime}}_0^G}}
\newcommand{\bl}{\mathrm{bl}}
\newcommand{\Knullbl}{\ensuremath{\Knull^{\bl}}}
\newcommand{\sm}{\mathrm{sm}}
\newcommand{\Knullsm}{\ensuremath{\Knull^{\sm}}}
\newcommand{\Knullschl}{\ensuremath{\widetilde{\operatorname{K}}_0}}
\newcommand{\Knullschlbl}{\ensuremath{\widetilde{\operatorname{K}}_0^{\bl}}}
\newcommand{\Knulldach}{\ensuremath{\widehat{\operatorname{K}_0}}}

\newcommand{\Vark}{\ensuremath{\mathit{Var}_k}}
\newcommand{\VarS}{\ensuremath{\mathit{Var}_S}}
\newcommand{\Var}[1]{\ensuremath{\mathit{Var}_{#1}}}
\newcommand{\Motk}{\ensuremath{\mathfrak{M}_k}}
\newcommand{\Bl}[2]{\ensuremath{\operatorname{Bl}_{#1}#2}}

\newcommand{\cN}{\ensuremath{\mathcal{N}}}
\newcommand{\cL}{\ensuremath{\mathcal{L}}}

\newcommand{\M}{\ensuremath{\mathcal{M}}}
\newcommand{\D}{\ensuremath{\mathcal{D}}}
\newcommand{\DG}{\ensuremath{\mathbb{D}}}
\newcommand{\aff}[1]{\ensuremath{\mathbb{A}^{#1}}}
\newcommand{\proj}[1]{\ensuremath{\mathbb{P}^{#1}}}
\newcommand{\Proj}{\ensuremath{\operatorname{\mathbb{P}}}}
\newcommand{\Spec}{\ensuremath{\operatorname{Spec}}}
\newcommand{\Hom}{\ensuremath{\operatorname{Hom}}}

\newtheorem{thm}{Theorem}[section]
\newtheorem{lem}[thm]{Lemma}
\newtheorem{prop}[thm]{Proposition}
\newtheorem{cor}[thm]{Corollary}

\newtheorem{listprop}[thm]{List of properties}
\theoremstyle{definition}
\newtheorem{defi}[thm]{Definition}
\theoremstyle{remark}
\newtheorem{rem}[thm]{Remark}

\newtheorem{step}{Step}

\begin{document}
\title[The universal Euler characteristic]{The universal Euler characteristic for varieties of
characteristic zero}
\author{Franziska Bittner}
\address{Mathematisch Instituut\\ Universiteit Utrecht\\
PO Box 80.010\\NL-3508 TA Utrecht\\Nederland}
\email{bittner@math.uu.nl}
\begin{abstract}
Using the weak factorization theorem 
we give a simple presentation for the value group of the
universal Euler characteristic with compact support for 
varieties of characteristic zero and describe the value group of the
universal Euler characteristic of pairs. 
This gives a new proof for the existence of natural Euler characteristics with values in
the Grothendieck group of Chow motives.
A generalization of the presentation to the relative setting allows us
to define duality and the six operations.
\end{abstract}
\maketitle

\section{Introduction}
Let $k$ be a field of characteristic zero. For an abelian group $A$,
an \emph{$A$-valued Euler 
characteristic with compact support}
for varieties over $k$ assigns to every such variety $X$ an element
$e(X)\in A$ only depending on the isomorphism
class of $X$ such that
$e(X) = e(X-Y)+e(Y)$ for $Y\subset X$ a closed subvariety. There is an
evident universal Euler characteristic for varieties. It takes values
in 
the \emph{naive Grothendieck group of varieties over $k$}, denoted by
$\Knull(\Vark)$,
which is the
free abelian group on isomorphism classes $[X]$ of varieties $X$ over
$k$ modulo the relations $[X]=[X-Y]+[Y]$ for $Y\subset X$ a closed
subvariety. It can be given a ring structure by taking products of
varieties. This group and its relative and equivariant analogues (to
be more precise, certain localizations and completions) are the
natural value groups for motivic integrals, as explained e.g. 
by Looijenga in~\cite{Looijenga} and by Denef and Loeser in~\cite{DenefLoeser}.

In this note we give a simple presentation of this
group in terms of smooth projective varieties. In fact we show that
$\Knull(\Vark)$ is the free abelian group on isomorphism classes of
smooth projective varieties modulo the relations $[\emptyset]=0$ and 
$[X]-[Y]=[\Bl{Y}{X}]-[E]$, where $Y\subset X$ is a
smooth closed subvariety, $\Bl{Y}{X}$ denotes the blow-up of $X$ along
$Y$ and $E$ denotes the exceptional divisor of this blow-up. 
The key ingredient is the factorization theorem proven by
W{\l}odarczyk in~\cite{Wlodarczyk} and Abramovich, Karu, Matsuki and 
W{\l}odarczyk in~\cite{Abramovich}: Any 
birational map between smooth complete irreducible varieties over $k$
can be factored into a sequence of blow-ups and blow-downs
with smooth centers.
It follows that
there is a unique Euler characteristic with compact support for $k$-varieties 
with values in the Grothendieck group of
Chow motives over $k$ which assigns to a smooth projective
variety the class of its Chow motive. This has already been proven 
by Gillet and Soul\'{e} in~\cite{GilletSoule} and by Guillen and 
Navarro Aznar in
\cite{GuillenNavarro}.

There is a similar, though somewhat more involved definition of an
Euler characteristic of pairs of varieties. We show that the
universal Euler characteristic for pairs has the same value ring
as the universal Euler characteristic with compact support, 
hence for every pair $(X,Y)$ we get an element 
$\langle X, Y \rangle$ in $\Knull(\Vark)$. In these terms $[X]$
corresponds 
to $\langle \overline{X} ,\overline{X} - X\rangle$, where 
$X\subset \overline{X}$ is a completion of $X$.

We are then in a position to show that there is a unique Euler
characteristic for pairs of $k$-varieties with values in the
Grothendieck group of Chow motives which for a smooth projective
variety $X$ assigns to the pair $(X,\emptyset)$ the Chow motive of $X$. 
This gives a new proof for a result by Guillen 
and Navarro Aznar in~\cite{GuillenNavarro}.

We also develop simple
presentations 
for the Grothendieck group of
varieties over a base variety since this
will make the theory much more effective. For instance, we are then
able to bring in the operations $f_*,f_!,f^*,f^!$, tensor product,
internal homomorphisms and duality familiar from Grothendieck's
duality theory
and to derive some expected formulae.

Furthermore we generalize the presentations to an equivariant setting.

I am indebted to Eduard Looijenga, my thesis advisor. I would like to
thank
Jaros{\l}aw W{\l}odarczyk for answering my questions about the weak
factorization theorem, and Jaros{\l}aw Wi{\'s}niewski and the
Mathematical Department of the University of Warsaw for their kind
hospitality during my stay as a guest of the Polish EAGER node.

\subsection*{Conventions:}
In the sequel $k$ denotes a field of characteristic zero. By a
\emph{variety} over $k$ we mean a reduced scheme of finite type over
$k$, not necessarily irreducible.

\section{The weak factorization theorem}
\label{fact}
The weak factorization theorem for varieties over a (not necessarily
algebraically closed) field $k$ of
characteristic zero is proven by W{\l}odarczyk in~\cite{Wlodarczyk} and
Abramovich, Karu, Matsuki and W{\l}odarczyk
in~\cite{Abramovich}. For convenience we recall parts of Theorem 0.3.1 
in~\cite{Abramovich} and 
the subsequent Remark 2.
\begin{thm}
Let $\phi:X_1 \dashrightarrow X_2$ be a birational map between
complete smooth connected varieties over $k$, let $U\subset X_1$ be an open
set where $\phi$ is an isomorphism. Then $\phi$ can be factored into a
sequence of blow-ups and blow-downs with smooth centers disjoint from
$U$:
There exists a sequence of birational maps 
$$
X_1 = V_0 \stackrel{\phi_1}{\dashrightarrow} V_1
 \stackrel{\phi_2}{\dashrightarrow} \dots
 \stackrel{\phi_i}{\dashrightarrow} V_i
 \stackrel{\phi_{i+1}}{\dashrightarrow} V_{i+1}
 \stackrel{\phi_{i+2}}{\dashrightarrow}
 \dots \stackrel{\phi_{l-1}}{\dashrightarrow}
 V_{l-1} \stackrel{\phi_l}{\dashrightarrow} V_l = X_2
$$
where $\phi = \phi_l \circ \phi_{l-1} \circ \cdots \phi_2 \circ
 \phi_1$, such that each factor $\phi_i$ is an isomorphism over $U$, and 
  $\phi_i:V_i \dashrightarrow V_{i+1}$ or $\phi_i^{-1}:V_{i+1}
\dashrightarrow V_{i}$ 
  is a morphism obtained by blowing up a smooth center
 disjoint from $U$ (here $U$ is identified with an open subset of $V_i$).

Moreover,
   there is an index $i_0$ such that for all $i\leq i_0$ the map 
   $V_i\dashrightarrow X_1$ is defined everywhere and projective, and for all 
   $i\geq i_0$ the map
   $V_i\dashrightarrow X_2$ is defined everywhere and projective.

If $X_1- U$ (respectively, $X_2- U$) is a
  simple normal crossings divisor, then the factorization can be chosen
 such that the inverse images of this
 divisor 
under $V_i\rightarrow X_1$
 (respectively, $V_i \rightarrow X_2$) are
 also  simple normal crossings divisors, 
and the centers of blowing up 
 have normal crossings
 with these divisors. 

If $\phi$ is equivariant under the action of a finite group,
 then the factorization can be chosen 
equivariantly.
\end{thm} 
We make some remarks on this theorem.
\begin{rem}
It is not explicitly stated in theorem 0.3.1 of~\cite{Abramovich}
that in the case of $X_1-U$ (respectively, $X_2-U$) a  simple normal
crossings divisor the inverse images of this divisor under
$V_i\rightarrow X_1$ (respectively, $V_i \rightarrow X_2$) are also
 simple normal crossings divisors, but it can be read off from the proof
(see 5.9 and 5.10).
\end{rem}

\begin{rem}
The completeness of the varieties $X_1$ and $X_2$ is not necessary as
has already been pointed out by Bonavero in~\cite{Bonavero}. It suffices that
$\phi$ is a \emph{proper} birational map $X_1
\dashrightarrow X_2$ between smooth connected varieties.  
This means that
the projections to $X_1$
and $X_2$ from the graph of $\phi$ are proper, which gives back usual
properness for an everywhere defined $\phi$. 
\end{rem}

\begin{rem}
We note that the theorem implies that if
$X_1$ and $X_2$ are varieties over a base variety $S$ and $\phi$ is
a map of $S$-varieties then the factorization is a factorization over
$S$. 
If $X_1$ and $X_2$ are projective over $S$ then so are the $V_i$'s.
\end{rem}

\begin{defi}
An action of a finite group on a variety is said to be \emph{good} if
every orbit is contained in an affine open set.
\end{defi}

\begin{rem}
Note that on a quasi-projective variety every action of a finite group
is good. More
generally, if $X\longrightarrow S$ is equivariant and $X$ is
quasi-projective over $S$, then the action on $X$ is good if the
action on $S$ is.
In particular, 
in case of an equivariant factorization with respect to a 
good action of a
finite group on $X_1$ and $X_2$ the action of on the $V_i$'s is automatically
good as they are projective over $X_1$, or $X_2$ respectively.
\end{rem}

\begin{rem}
\label{trans}
If $X$ is a not necessarily connected variety with a (good)
$G$-action such that $G$ acts transitively on the connected
components of $X$ one can identify $X$ with the induced variety 
$G\times^{G_0} X_0= (G \times X_0) / G_0$, where
$X_0$ is a connected component of $X$ and $G_0$ is the stabilizer of
$X_0$. Here the action of $G_0$ is given by $h(g,x)=(gh^{-1}, hx)$.

Using this and the equivariant version of the factorization theorem
for connected smooth varieties also yields an equivariant factorization for
smooth, not necessarily connected varieties where $G$ acts transitively on
the connected components.
\end{rem}

\section{The universal Euler characteristic with compact support}
\label{univkomp}
Let $\Knull(\Vark)$ be the free abelian group on isomorphism classes
$[X]$ of
varieties $X$ over $k$ where we impose the relation
$[X]=[X-Y]+[Y]$ for $Y\subset X$
a closed subvariety (so in particular $[\emptyset]=0$).  
This group is also called \emph{the
(naive) Grothendieck group of $k$-varieties}. It is the value group of
the universal Euler characteristic with compact support.
It is a commutative ring for the multiplication induced by 
the product of
$k$-varieties, because 
the subgroup divided out actually is an ideal with respect to this
multiplication. 

\begin{thm}
\label{preskomp}
The Grothendieck group of $k$-varieties has the following alternative
presentations: 
\begin{itemize}
\item[(sm)] As the abelian group generated by the isomorphism classes
of smooth varieties over $k$ subject to the relations $[X]=[X-Y]+[Y]$, where
$X$ is smooth and $Y\subset X$ is a smooth closed subvariety,
\item[(bl)] as the abelian group generated by the isomorphism classes
of smooth complete $k$-varieties subject to the relations
$[\emptyset]=0$ and $[\Bl{Y}{X}] -[E] = [X] -[Y]$, where 
$X$ is smooth and complete, $Y\subset X$ is a closed smooth
subvariety, $\Bl{Y}{X}$ is the blow-up of $X$ along $Y$ and $E$ is the
exceptional divisor of this blow-up.
\end{itemize}
Moreover, we get the same group if in case (sm) we restrict to
quasi-projective varieties or if in case (bl) we restrict to projective 
varieties. We can also restrict to connected varieties in both
presentations.  
\end{thm}
\begin{rem}
The subgroups divided out are again ideals with respect to the product
of $k$-varieties.
\end{rem}

\begin{defi}
Let $\tate=[\aff{1}]$ denote the class of the affine line.
Define the \emph{naive motivic ring of $k$-varieties}, $\M_k$, as the
localization $\Knull(\Vark)[\tate^{-1}]$.
\end{defi}

\begin{cor}
\label{dual}
We have a ring involution $\D_k$ of $\M_k$ that sends $\tate$ to
$\tate^{-1}$ and is characterized by the property that it sends the
class of a complete connected smooth variety $X$ to $\tate^{-\dim X}[X]$.
\end{cor}
We call this involution the \emph{duality map}.

We first deduce Corollary~\ref{dual} from Theorem~\ref{preskomp}.

A vector bundle $V\longrightarrow X$ of rank $n$ on 
a variety $X$ is locally trivial by definition and hence 
its class $[V]$ in $\Knull(\Vark)$ is equal to $\tate^n [X]$. 
Similarly the
class of its projectivization $[\Proj(V)]$ equals $[\proj{n-1}][X]=
(1+\tate + \dots + \tate^{n-1})[X]$.

\begin{lem}
\label{fund2}
Let $X$ be a smooth connected variety and $Y\subset X$ a smooth
connected subvariety, let $d$ denote the codimension of $Y$ in $X$. 
Let $E$ be the exceptional divisor of the blow-up $\Bl{Y}{X}$ of $X$
along $Y$. Then
$[\Bl{Y}{X}]-\tate[E] = [X]-\tate^d[Y]$ in $\Knull(\Vark)$. 
\end{lem}
\begin{proof}
In $\Knull(\Vark)$ we have the fundamental relation
$[\Bl{Y}{X}] - [E] = [X]-[Y]$. 
Furthermore $[E] = (1+ \tate + \dots + \tate^{d-1})[Y]$, thus $(1-\tate)[E] =
(1-\tate^d)[Y]$.
Adding this to the fundamental relation finishes the proof. 
\end{proof}

\begin{proof}[Proof of Corollary~\ref{dual}]
Using the presentation (bl) and Lemma~\ref{fund2} we can define 
a group homomorphism $\Knull(\Vark)\longrightarrow \M_k$ by sending 
the class of a smooth connected complete variety $X$ to $\tate^{-\dim X}[X]$.
This morphism is multiplicative and maps $\tate$ to $\tate^{-1}$,
hence it can be extended uniquely to a ring endomorphism 
$\D_k$ of $\M_k$. Obviously $\D_k\D_k = \id{\M_k}$.
\end{proof}

We now deduce the presentation (sm) of Theorem~\ref{preskomp}.
For this purpose let us provisionally introduce the group
$\Knullsm(\Vark)$, defined as the free abelian group on isomorphism classes
$[X]_{\sm}$ of smooth varieties modulo the relations for smooth 
closed subvarieties $Y$ of smooth varieties $X$.
It is a commutative ring with respect to the product of
varieties.

\begin{step}
The ring homomorphism $\Knullsm(\Vark) \longrightarrow \Knull(\Vark)$, 
$[X]_{\sm} \mapsto [X]$
is an isomorphism.
\end{step} 

\begin{proof} To construct an inverse for any variety $X$ we stratify
$X=\bigsqcup_{N\in\cN} N$ such that $N$
smooth 
and equidimensional
and $\overline{N}$ a union of strata for all $N\in\cN$. 
Consider the expression $\sum_{N\in\cN} [N]_{\sm}$ in $\Knullsm(\Vark)$.
If $X$ is smooth itself 
$\sum_{N\in\cN} [N]_{\sm}$ equals $[X]_{\sm}$ as can be seen by
induction on the number of elements of $\cN$: Let $N\in\cN$ be an
element of minimal dimension, then $[X]_{\sm}=[X-N]_{\sm} + 
[N]_{\sm}$, and by the induction hypothesis 
$[X-N]_{\sm}= \sum_{N'\in\cN-\{N\}} [N']_{\sm}$.

For two stratifications $\cN$ and $\cN'$ of $X$ we can always
find a common refinement $\cL$. The above argument shows that 
for $N\in\cN$ we get $\sum_{N\supset L\in\cL} [L]_{\sm} = [N]_{\sm}$.
Hence
$\sum_{L\in\cL} [L]_{\sm}$ is equal to $\sum_{N\in\cN}[N]_{\sm}$ and 
analogously it equals $\sum_{N\in\cN'}[N]_{\sm}$,
therefore $\sum_{N\in\cN} [N]_{\sm}$ is independent of the
choice of the stratification. Thus we can set 
$e(X):=\sum_{N\in\cN} [N]_{\sm}$.

If $Y\subset X$ is a closed subvariety we can find a stratification
for which $Y$ is a union of strata which yields
$e(X)=e(X-Y) + e(Y)$, hence $e$ is an Euler characteristic with
compact support and factors through
$\Knull(\Vark)$. The induced map on 
$\Knull(\Vark)$ 
obviously is an inverse for $\Knullsm(\Vark)\longrightarrow \Knull(\Vark)$.
\end{proof}

Decomposing into connected components and noting that instead of
cutting a smooth closed subvariety $Y$ out of a smooth connected
variety $X$ we can also take out the connected components of $Y$
one by one shows that we can restrict to smooth connected varieties in
the presentation (sm).

Stratifying by smooth quasi-projective varieties shows that we can
restrict to smooth (connected) quasi-projective varieties.

In the rest of this section we establish the presentation (bl) of
Theorem~\ref{preskomp}.

Again we introduce an auxiliary group:
Let $\Knullbl(\Vark)$ be the free abelian group on isomorphism classes
$[X]_{\bl}$ of smooth complete varieties $X$ over $k$ 
modulo the relations for blow-ups of smooth complete varieties $X$
along smooth closed subvarieties $Y$ and the relation 
$[\emptyset]_{\bl}=0$
(then $[X\sqcup
Y]_{\bl}=[X]_{\bl} + [Y]_{\bl}$, which can be seen by blowing up along
$Y$).

Decomposing into connected components and noting
that the blow-up along a disjoint union is the successive blow-up along
the connected components one sees that this can also be described as
the free abelian group on isomorphism classes $[X]_{\bl}$ of connected
smooth complete varieties with imposed relations
$[\emptyset]_{\bl}=0$ and $[\Bl{Y}{X}]_{\bl} -[E]_{\bl} = [X]_{\bl}
-[Y]_{\bl}$, where $Y\subset X$ is a connected closed smooth
subvariety. 

Also $\Knullbl(\Vark)$ carries a commutative ring structure induced by 
the product of varieties.

\begin{step}
The ring homomorphism $\Knullbl(\Vark) \longrightarrow \Knull(\Vark)$ which 
sends $[X]_{\bl}$ to $[X]$ is an isomorphism.
\end{step} 

\begin{proof}
Again we construct an inverse. Using the presentation (sm) in
Theorem~\ref{preskomp} we see that defining an Euler characteristic
with compact support is equivalent to defining an Euler characteristic
with compact support for smooth connected varieties.

Let $X$ be a smooth connected variety, let $X\subset
\overline{X}$ be a smooth completion with $D=\overline{X}-X$ a simple normal 
crossings divisor. Let
$D^{(l)}$ be the normalization of the $l$-fold intersections of $D$,
where $D^{(0)}$ is understood to be $\overline{X}$ (so $D^{(l)}$ is the
disjoint union of the $l$-fold intersections of the irreducible
components of $D$). 
Consider the expression $\sum (-1)^l [D^{(l)}]_{\bl}$ in $\Knullbl(\Vark)$.

We first show that this expression is independent of the choice of the
completion:

Let $X\subset \overline{X}'$ and $X \subset \overline{X}$ be two
smooth completions of $X$ with $\overline{X}-X=D$ and
$\overline{X}'-X=D'$ simple normal crossings divisors. Due to the weak factorization 
theorem the birational map 
$\overline{X}'\dashrightarrow\overline{X}$ can be factored into a
sequence of blow-ups and blow-downs 
with smooth centers disjoint from $X$ which have normal crossings with
the complement of $X$. Hence we may assume that
$\overline{X}' = \Bl{Z}{\overline{X}}$ with $Z\subset D$ smooth and
connected such that $Z$ has normal crossings with
$D$.

Let $D_0$ be the irreducible component of $D$ containing $Z$ and let
$\{D_i\}_{i\in I}$ be the remaining irreducible components.
Then the irreducible components of $D'$ are $D'_i=\Bl{Z\cap D_i}D_i$
(where $i\in\{0\}\cup I$) and the exceptional divisor $E$ of the
blow-up.
For $K\subset \{0\}\cup I$ we put $D_K:=\bigcap_{j\in K} D_j$ (where
$D_{\emptyset}$ is understood to be $\overline{X}$), 
$D'_K:=\bigcap_{j\in K} D'_j$, $Z_K:=Z\cap D_K$ and $E_K:=E\cap D'_K$.
As $Z$ has normal crossings with $D$ we get $D'_K=\Bl{Z_K}{D_K}$ with
exceptional divisor $E_K$, hence we have
$$[D'_K]_{\bl}-[E_K]_{\bl}= [D_K]_{\bl}-[Z_K]_{\bl}.$$
Denote by 
$E^{(l)}$ the preimage of $E$ in  ${D'}^{(l)}$ and by $Z^{(l)}$
the preimage of $Z$ in $D^{(l)}$. Then 
for 
$l=0,\dots,n$ the preceding identity yields
\begin{align*}
[{D'}^{(l)}]_{\bl}&=\sum_{|K|=l} [D'_K]_{\bl} +
\sum_{|K|=l-1}[E_K]_{\bl}\\
&= \sum_{|K|=l} ([D_K]_{\bl}+[E_K]_{\bl} - [Z_K]_{\bl}) +
\sum_{|K|=l-1}[E_K]_{\bl}\\
&= [D^{(l)}]_{\bl} + [E^{(l)}]_{\bl} - [Z^{(l)}]_{\bl} +
[E^{(l-1)}]_{\bl}, 
\end{align*} 
(for $l=0$ the last term is zero).
As $Z\subset D_0$
we get $Z_{\{0\}\cup K} = Z_K$ for $K\subset I$, thus 
$\sum (-1)^l[Z^{(l)}]_{\bl} = 0$.
Taking the alternating sum hence yields
$$\sum (-1)^l [{D'}^{(l)}]_{\bl} = \sum (-1)^l [D^{(l)}]_{\bl}.$$

Therefore we can set $e(X) := \sum (-1)^l [D^{(l)}]_{\bl}$.

We have to check that $e(X)=e(X-Y)+ e(Y)$ for $Y\subset X$ a 
connected closed smooth subvariety of a connected smooth variety $X$. 
We choose $\overline{X}\supset X$
smooth and complete such that $D = \overline{X} - X$ is a simple normal
crossings divisor and such that the closure $\overline{Y}$ of $Y$ in
$\overline{X}$ is also smooth and has normal crossings with $D$
(we can take first a smooth completion of $X$ with boundary a
simple normal crossings divisor and then an embedded resolution of
the closure of $Y$ compatible with this divisor -- compare
e.g. Section~1.2 of~\cite{Abramovich}). 
In particular $D\cap \overline{Y}$
is a simple normal crossings divisor in $\overline{Y}$. Denote the
irreducible components of $D$ by $\{D_i\}_{i\in I}$, for $K\subset I$
let $D_K$ be defined as above, let $Y_K:=\overline{Y}\cap D_K$ and
$Y^{(l)}=\bigsqcup_{|K|=l} Y_K$. Then $e(Y) = \sum
(-1)^l[Y^{(l)}]_{\bl}$.

Let $\widetilde{X}:= \Bl{\overline{Y}}{\overline{X}}$ and denote the
exceptional divisor by $E$.
Denote the proper transform of $D_i$ by $\widetilde{D}_i$.
The complement $\widetilde{D}$ of $X-Y$ in $\widetilde{X}$ 
is the simple normal crossings divisor $\bigcup \widetilde{D}_i \cup E$.
If $E_K:=E\cap \widetilde{D}_K$
then as above $\widetilde{D}_K$ is the blow-up of $D_K$ along $Y_K$
with exceptional divisor $E_K$ and hence 
$$[\widetilde{D}_K]_{\bl}-[E_K]_{\bl} = [D_K]_{\bl} - [Y_K]_{\bl}.$$
Thus
\begin{align*}
[\widetilde{D}^{(l)}]_{\bl} &= \sum_{|K|=l} [\widetilde{D}_K]_{\bl} +
\sum_{|K|=l-1} [E_K]_{\bl}\\
&=\sum_{|K|=l} ([D_K]_{\bl} -[Y_K]_{\bl} + [E_K]_{\bl}) 
+ \sum_{|K|=l-1} [E_K]_{\bl}\\
&= [D^{(l)}]_{\bl} - [Y^{(l)}]_{\bl} + [E^{(l)}]_{\bl} +
[E^{(l-1)}]_{\bl}, 
\end{align*}
where $E^{(l)}$ denotes the preimage of $E$ in $\widetilde{D}^{(l)}$.
Taking the alternating sum yields $e(X-Y)=e(X)-e(Y)$.
Hence $e$ induces a morphism $\Knull(\Vark)\longrightarrow
\Knullbl(\Vark)$ which clearly is an inverse for the mapping
$\Knullbl(\Vark)\longrightarrow \Knull(\Vark)$.
\end{proof}

Using the fact that we can restrict to quasi-projective generators in 
the presentation (sm) of Theorem~\ref{preskomp} and that a connected smooth
quasi-projective variety has a smooth projective
simple normal crossings completion we see that we
can restrict to projective generators in the description (bl) of
Theorem~\ref{preskomp}.

\section{The universal Euler characteristic of pairs}
\label{univpaar}
Let $A$ be an abelian group. Then an \emph{$A$-valued Euler
characteristic} is a mapping which associates to each pair $(X,Y)$ of
varieties over $k$ with $Y\subset X$ a closed subvariety an element
$\chi(X,Y)\in A$ (where we denote $\chi(X,\emptyset)$ by $\chi(X)$)
such that the following properties hold:
\begin{itemize}
\item \emph{Excision}: If $f:X' \longrightarrow X$ is proper and
$Y\subset X$ is a closed subvariety such that $f$ induces an isomorphism
$X'-f^{-1}Y\cong X-Y$, 
then $\chi(X',f^{-1}Y) = \chi (X,Y)$.
\item \emph{Gysin}: If $X$ is smooth and connected and $D\subset X$ is a smooth
divisor, then $\chi (X- D) = \chi (X) - \chi (\proj{1}\times D,
\{\infty\}\times D)$.
\item \emph{Exactness}: If $X\supset_{\text{closed}} Y
\supset_{\text{closed}} X$, then $\chi(X,Z) = \chi(X,Y) + \chi(Y,Z)$.
\end{itemize}
Excision implies that
$\chi(X,Y)$ only depends on the isomorphism class of the pair $(X,Y)$.
Using exactness we get $\chi(X,Y) = \chi(X) - \chi(Y)$.
Exactness and excision yield $\chi(X\sqcup Y)= \chi(X\sqcup
Y,Y)+\chi(Y) = \chi(X)+ \chi(Y)$ and $\chi(\emptyset)=0$.
If $X\subset W$ is an open embedding with $W$ complete, then the
the excision property implies that
$\chi(W, W-X)$ is independent of the choice of
the open embedding. 
We denote $\chi(W,W-X)$ by $\chic(X)$.
For $Y\subset X$ closed we have
\begin{multline*}
\chic(X) = \chi(W, W-X)\\ = \chi(W,
\overline{Y}\cup (W-X)) + \chi(\overline{Y} \cup
(W-X), W-X) = \chic(X-Y) + \chic(Y),
\end{multline*}
where by $\overline{Y}$ we denote the closure of $Y$ in $W$.
Hence $\chic$ factors through $\Knull(\Vark)$.

\begin{defi}
The 
\emph{value ring of the universal Euler characteristic for pairs
of $k$-varieties} is defined to be the free abelian group on isomorphism
classes of
pairs $(X,Y)$ with $Y\subset X$ closed modulo the relations for an
Euler characteristic. 
We denote it by $\Knullschl(\Vark)$. The class of a pair $(X,Y)$ in
$\Knullschl(\Vark)$ is denoted by
$\langle X, Y \rangle$, where we write $\langle X \rangle$ shorthand for
$\langle X, \emptyset\rangle$.
\end{defi}
The above construction yields a group homomorphism
$$\chic:\Knull(\Vark) \longrightarrow \Knullschl(\Vark)$$

Indeed $\Knullschl(\Vark)$ carries a ring structure by defining the
product of 
$\langle X,Y\rangle$ and $\langle X', Y'\rangle$ to be 
$\langle X\times X', X\times Y'\cup Y\times X'\rangle$:

To see that this is compatible with the Gysin relation we introduce an
auxiliary group
$\Knulldach(\Vark)$, defined as the free 
abelian group on isomorphism classes of pairs of varieties over $k$ 
modulo excision and exactness. We denote the class of a pair $(X,Y)$
in $\Knulldach(\Vark)$ by $\{X,Y\}$, where we write
$\{X\}$ for $\{X,\emptyset\}$. 
We get a ring structure on $\Knulldach(\Vark)$ by setting
$\{ X,Y\} \{ X',Y'\} 
=\{ X\times X', X\times Y'\cup Y\times X'\}$. Noting that
$\{X,Y\}=\{X\} - \{Y\}$ and 
$\{X,S\} = \{\widetilde{X}, \widetilde{S}\}$ for an arbitrary $X$
where $S$ is the singular locus of $X$ and $\widetilde{X}$ a
resolution of singularities with $\widetilde{S}$ the inverse image of
$S$ we get $\{X\}=\{\widetilde{X}\}-\{\widetilde{S}\}+\{S\}$, where the
dimension of $\widetilde{S}$ and $S$ is strictly less than the
dimension of $X$. Proceeding inductively on the dimension of $X$ we
can write every class $\{X,Y\}$ as a linear combination of classes of
smooth varieties.
Hence the subgroup generated by expressions of the form 
$\{X- D\} - \{X\} + \{\proj{1}\times D, \{\infty\}\times D\}$ with $X$
smooth and $D$ a smooth divisor on $X$ actually is an ideal.

Using the multiplicative structure of $\Knullschl(\Vark)$ the Gysin
relation can be rewritten as 
$$\langle X-D\rangle = \langle X\rangle -\langle \proj{1},
\{\infty\}\rangle\langle D\rangle.$$

The group homomorphism $\chic$ defined above actually is 
a ring homomorphism.

\begin{thm}
\label{prespaar}
The map $\chic:\Knull(\Vark) \longrightarrow \Knullschl(\Vark)$
is a ring isomorphism. If $X$ is a smooth connected $k$-variety and
$(\overline{X},D)$ is a simple normal crossings completion of $X$ over $k$,
then the inverse of $\chic$ assigns to $\langle X\rangle$ the element 
$\sum (-\tate)^l [D^{(l)}]$. Its image in $\M_k$ is also equal to
$\D_k(\tate^{-\dim X}[X]) = \tate^{\dim X}\D_k([X])$.
\end{thm}

Before we prove this we give an application of the 
theorems~\ref{preskomp} and~\ref{prespaar}.

Let $\Motk$ denote the category of Chow $k$-motives
(see e.g. the introduction~\cite{Scholl} by
Scholl for the definition of $\Motk$).
For
a smooth connected projective variety $X$ over $k$, denote by $h(X)$ 
the motive of $X$ and its class in
$\Knull(\Motk)$ by $[h(X)]$.  
Denote the class $[h(\proj{1})] - [h(\Spec k)]$
of the Tate motive by
$\tatemot$, and by $[h(X)]^\vee = \tatemot^{-\dim X}\otimes [h(X)]$ the
class of the dual motive of $h(X)$.

\begin{cor}
\label{chow}
There is a ring homomorphism (in fact the unique
group homomorphism)
$\chicmot:\Knull(\Vark) \longrightarrow \Knull(\Motk)$
which sends the class of a smooth projective variety $X$ to $[h(X)]$. 
The class of the affine line is mapped to the class of the Tate motive 
by this morphism. 
Likewise, there is a ring homomorphism (in fact the unique group homomorphism)
$\chimot:\Knullschl(\Vark) \longrightarrow \Knull(\Motk)$
which sends the class of a smooth projective variety $X$ to $[h(X)]$.
The two are related by the property that
for $X$ smooth and connected we have 
$(\chimot\langle X\rangle)^\vee= \tatemot^{-\dim X}\otimes\chicmot[X]$.
\end{cor}

\begin{proof}[Proof of Corollary~\ref{chow}]
For $Y\subset X$ a smooth closed subvariety of a smooth projective
variety we have $[h(\Bl{Y}{X})]-[h(E)] = [h(X)]-[h(Y)]$, where $E$
denotes the exceptional divisor. We now use the presentation (bl) of
Theorem~\ref{preskomp} to see that there is a unique group homomorphism
$\Knullschl(\Vark) \longrightarrow \Knull(\Motk)$
which sends the class of a smooth projective variety $X$ to $[h(X)]$
and therefore is a ring homomorphism. As $\tate$ is mapped to $\tatemot$
it can be extended to $\M_k$.

Furthermore due to Theorem~\ref{prespaar} 
there is a unique group homomorphism $\chimot$ (in fact a
ring homomorphism) which makes
$$\xymatrix@R=1ex{\Knull(\Vark)\ar[rd]^{\chicmot}\ar[dd]_{\chic}&\\
            &\Knull(\Motk)\\
          \Knull(\Vark \ar[ru]^{\chimot})&}$$
commutative.

Now let $X$ be a smooth connected variety.
As the image of the inverse of $\langle X \rangle$ under $\chic$ in
$\M_k$ is $\D_k(\tate^{-\dim X} [X])$
and as
$$\xymatrix{\M_k \ar[r]^-{\chicmot} \ar[d]_{\D_k}& \Knull(\Motk)\ar[d]^{\vee}\\
             \M_k \ar[r]^-{\chicmot}&\Knull(\Motk)}$$
is commutative we conclude 
that 
$(\chimot\langle X\rangle)^\vee= \tatemot^{-\dim X}\otimes\chicmot[X]$.
\end{proof}

The rest of this section is devoted to the proof of
Theorem~\ref{prespaar}.

For this purpose
we give a more convenient presentation of
$\Knullschl(\Vark)$. 

For the moment let  $\Knullschlbl(\Vark)$ be 
the free abelian group on isomorphism classes of pairs of 
\emph{smooth} varieties (denoted by
$\langle X,Y\rangle_{\bl}$, where $\langle X \rangle_{\bl}$  stands for
$\langle X, \emptyset\rangle_{\bl}$), with the same relations as for
$\Knullschl(\Vark)$ except that we impose 
the excision
relation only for blow-ups at smooth centers 
(as $\langle X, \emptyset\rangle_{\bl} = \langle X\sqcup Y,
Y\rangle_{\bl}$ by blowing up $X\sqcup Y$ along $Y$ we get by exactness
$\langle X\sqcup Y\rangle_{\bl} = \langle X\rangle_{\bl} + \langle
Y\rangle_{\bl}$ and $\langle \emptyset \rangle_{\bl} = 0$).

\begin{prop}
\label{presbl}
The obvious group homomorphism $\Knullschlbl(\Vark) \longrightarrow 
\Knullschl(\Vark)$ is an
isomorphism. In other words, 
the value group $\Knullschl(\Vark)$ of the universal Euler
characteristic for pairs of $k$-varieties
has a presentation as the free group on isomorphism classes of pairs
of 
smooth $k$-varieties, modulo 
exactness, Gysin and excision for blow-ups at smooth centers.
\end{prop}
The strategy for the proof of the proposition is
to construct an inverse by defining a map 
$\Psi:\Knulldach(\Vark)\longrightarrow \Knullschlbl(\Vark)$  
that is the obvious map on a smooth pair $\{X,Y\}$ and
factors through the Gysin relation.
It is enough to define $\Psi\{X\}$ for all varieties $\{X\}$ and then to set
$\Psi\{X,Y\}= \Psi\{X\}-\Psi\{Y\}$. Then $\Psi$ will automatically
fulfill exactness.
For a smooth variety $X$ we want $\Psi\{X\}$ to be equal to 
$\langle X\rangle_{\bl}$.

We proceed by induction on the dimension. To be more precise
we use that $\varinjlim \Knulldach(\Vark)_n\cong\Knulldach(\Vark)$, 
where $\Knulldach(\Vark)_n$ denotes the free
abelian group on isomorphism classes of pairs of varieties of
dimension at most $n$ modulo excision and exactness.
In dimension zero we just set $\Psi\{X\} = \langle X\rangle_{\bl}$.
Suppose now that $\Psi$ has already been defined for pairs $\{X,Y\}$
with $\dim X < n$, that it factors through the excision relation for such
pairs (in particular it is additive on disjoint unions in dimensions
smaller than $n$) and that $\Psi\{X\}=\langle X\rangle_{\bl}$ for $X$ 
smooth and $\dim X< n$.
In the following five steps we extend $\Psi$ to $\Knulldach(\Vark)_n$.
\setcounter{step}{0}
\begin{step} 
\label{Ausschneidungglatt}
Suppose $f:X'\longrightarrow X$ is proper with $Y'=f^{-1}Y$, suppose
that $f$ induces an isomorphism $X'-Y'\cong X-Y$.  Suppose $\dim X,
\dim X'\le n$ and $\dim Y,\dim Y'< n$, $X$ and $X'$ smooth.  Then
$\langle X'\rangle_{\bl}-\Psi\{Y'\}=
\langle X\rangle_{\bl}-\Psi\{Y\}$.
\end{step}
\begin{proof}
Suppose first that $X$ and $X'$ are connected.
Due to the factorization theorem $f$ can be
factored into a sequence of blow-ups and blow-downs with smooth
centers disjoint from $X-Y$. So we can assume that $X'\longrightarrow X$
is a blow-up along
$Z\subset X$ a smooth closed subvariety with $Z\subset Y$. Let
$E\subset Y'$ be
the exceptional divisor of this blow-up. We have
$\langle X'\rangle_{\bl} -\langle E\rangle_{\bl} = 
\langle X\rangle_{\bl} - \langle Z\rangle_{\bl}$
and furthermore $\Psi\{Y'\}-\langle E\rangle_{\bl} = \Psi\{Y\} -
\langle Z\rangle_{\bl}$ by
induction hypothesis. Hence $\langle X'\rangle_{\bl}-\Psi\{Y'\}=\langle
X\rangle_{\bl} -\Psi\{Y\}$.

Now in general $X$ can be written as the disjoint union of connected
components $X_1 \sqcup \dots \sqcup X_l \sqcup \dots \sqcup X_n$ such
that $Y_i:= X_i\cap Y$ is a proper subset of $X_i$ for $i\le l$ and
$Y_i = X_i$ for $i>l$. In particular $Y_i$ is smooth for $i>l$, and as 
$\dim Y_i < n$ the induction hypothesis yields $\Psi\{Y_i\} = \langle
Y_i\rangle_{\bl} = \langle X_i\rangle_{\bl}$ and hence 
$$\langle X \rangle_{\bl} - \psi\{Y\}=  
\sum_{i=1}^l\langle X_i \rangle_{\bl} - \psi\{Y_i\}.$$
Decomposing $X'= X'_1\sqcup \dots \sqcup X'_{l'} \sqcup \dots \sqcup
X'_{n'}$ in the same way and setting $Y'_i=X'_i\cap Y$ yields
$$\langle X' \rangle_{\bl} - \psi\{Y'\}=  
\sum_{i=1}^{l'}\langle X'_i \rangle_{\bl} - \psi\{Y'_i\}.$$
As $f$ induces an isomorphism 
$\bigsqcup_{i=1}^{l'} X'_i - Y'_i 
\cong \bigsqcup_{i=1}^l X_i - Y_i$ it follows that $l=l'$ and we may
assume that $f$ induces $(X'_i,Y'_i)\longrightarrow (X_i,Y_i)$ for
$i\le l$, hence the two sums are equal.  
\end{proof}

Now let $X$ be an arbitrary $n$-dimensional variety. 
Choose $\pi:\widetilde{X}\longrightarrow X$ proper with
$\widetilde{X}$ smooth such that $\pi$
induces an isomorphism over the smooth locus of $X$. Let $S\subset X$
be the singular locus of $X$, let $\widetilde{S}\subset \widetilde{X}$
be the inverse image of $S$.
As we want $\Psi\{\widetilde{X},\widetilde{S}\}$ to be equal to
$\Psi\{X,S\}$ we set
$\Psi\{X\}:=\langle\widetilde{X}\rangle_{\bl}-\Psi\{\widetilde{S}\}+
\Psi\{S\}$ (the right hand
side has already been defined as $\widetilde{X}$ is smooth and $\dim S,
\dim \widetilde{S} < n$).

Since any two choices of $\widetilde{X}\longrightarrow X$ are
dominated by a third one (take for example the closure of $X-S$ which is
diagonally embedded in the product and resolve the singularities),
Step~\ref{Ausschneidungglatt} implies that 
$\langle\widetilde{X}\rangle_{\bl}-\Psi\{\widetilde{S}\}+
\Psi\{S\}$  is independent
of this choice.
It is clear that $\Psi$ is additive with respect to disjoint unions
and that for smooth $X$ we get $\Psi\{X\}=\langle X\rangle_{\bl}$. 

We have to prove that $\Psi$ factors through excision for pairs of
varieties of dimension $\le n$.

\begin{step}
The map $\Psi$ factors through excision for 
$f:(X',Y')\longrightarrow (X,Y)$ 
with $X$ and $X'$ both smooth of dimension $\le n$: 
\end{step}

\begin{proof}
Like in the proof of Step~\ref{Ausschneidungglatt} we can decompose
$X$ and $X'$ into connected components. Now use $\Psi\{X_i\} = 
\langle X_i \rangle_{\bl}$ and $\Psi\{X'_i\} = 
\langle X'_i \rangle_{\bl}$ and apply 
Step~\ref{Ausschneidungglatt} to the $X_i$ and $X'_i$ with $i\le l$.
\end{proof}

\begin{step}
The map $\Psi$ factors through excision for 
$f:(X',Y')\longrightarrow (X,Y)$ with $\dim X , \dim
X' \le n$, $X'$ smooth and $\dim Y  < n$.
\end{step}

\begin{proof}
Let $S\subset Y$ be the singular locus of $X$, let $\widetilde{X}
\longrightarrow X$ be a resolution of singularities, let
$\widetilde{Y}$ be the inverse image of $Y$ and $\widetilde{S}$ the
inverse image of $S$. Then by definition $\Psi\{X,S\} = \Psi
\{\widetilde{X},\widetilde{S}\}$ and by induction $\Psi\{Y,S\} =
\Psi\{\widetilde{Y}, \widetilde{S}\}$. Hence $\Psi\{X,Y\} = \Psi
\{\widetilde{X}, \widetilde{Y}\}$. 

Let $\widehat{X}$ be a resolution of
singularities of the closure of $X-Y$ in 
$\widetilde{X}\times X'$ (we embed $X-Y$ diagonally). 
We have a commutative diagram 
$$\xymatrix{\widehat{X} \ar[r] \ar[d] &
X'\ar[d] \\ \widetilde{X} \ar[r] & X}$$ 
of proper birational maps which
yield isomorphisms over $X-Y$.  
Let $\widehat{Y}$ be the inverse image of $Y$ in $\widehat{X}$. 
Then $\Psi\{\widehat{X},\widehat{Y}\} =
\Psi\{X',Y'\}$ and $\Psi\{\widehat{X},\widehat{Y}\}= \Psi
\{\widetilde{X},\widetilde{Y}\}$ as $\widehat{X}$, $\widetilde{X}$ and
$X'$ are smooth.
\end{proof}

\begin{step}
The map $\Psi$ factors through excision for $f:(X_1,Y_1)\longrightarrow
(X_2,Y_2)$ with $\dim X_1, \dim X_2 \le n$ and $\dim Y_1, \dim Y_2 <
n$.
\end{step}

\begin{proof}
For $i=1,2$ 
let $S_i$ be the singular locus of $X_i$, let $Z_i = Y_i\cup S_i$.
Then $f^{-1}(Z_2) = Z_1$ and $f$ induces an isomorphism of smooth varieties
$X_1-Z_1\cong X_2-Z_2$ (as the singular locus of $X_i-Y_i$ is
$S_i-Y_i$).
Let $\widetilde{X_i}$ be a resolution of singularities of $X_i$ 
that is an isomorphism over $X_i-Z_i$: so if $\widetilde{Z_i}$ 
is the preimage of $Z_i$ in $\widetilde{X_i}$ then
$\widetilde{X_i}-\widetilde{Z_i} \cong X_i - Z_i$. 
Let $X_3$ be the closure of $X_1-Z_1$ (embedded diagonally) in
$\widetilde{X_1}\times\widetilde{X_2}$. Let $\widetilde{X_3}$ a
resolution of singularities of $X_3$ (that is an isomorphism over
$X_1-Z_1$). 
We thus have a commutative diagram
$$\xymatrix@C=1ex{&\widetilde{X_3}\ar[dl]_{\pi_1} \ar[dr]^{\pi_2}&\\
\widetilde{X_1}\ar[d]&&\widetilde{X_2}\ar[d]\\
X_1\ar[rr]&& X_2}$$ 
where $\pi_1^{-1}(\widetilde{Z_1})=\pi_2^{-1}(\widetilde{Z_2})=:\widetilde{Z_3}$ and
$\pi_1$, $\pi_2$ induce isomorphisms outside $\widetilde{Z_3}$.
As $\widetilde{X_i}$ are smooth for $1\le i\le 3$ and $\dim Z_1, \dim
Z_2 < n$ we have 
$$\Psi\{X_1,Z_1\}
=\Psi\{\widetilde{X_1},\widetilde{Z_1}\} 
= \Psi\{\widetilde{X_3},\widetilde{Z_3}\}\\
=\Psi\{\widetilde{X_2},\widetilde{Z_2}\}
=\Psi\{X_2,Z_2\}.$$
The induction hypothesis gives 
$$\Psi\{Z_1,Y_1\}=\Psi\{Z_2,Y_2\}.$$
Adding these two equations up yields $\Psi\{X_1,Y_1\}=\Psi\{X_2,Y_2\}.$
\end{proof}

\begin{step}
The map $\Psi$ factors through excision for $f:(X',Y')\longrightarrow
(X,Y)$ with $\dim X , \dim X' \le n$ .
\end{step}

\begin{proof}
Let $X_1,\dots,X_l$ be the irreducible components of $X$, 
let $Y_i := X_i \cap Y$. Let $B:=
\bigcup_{Y_i \ne X_i} X_i$, let $N:=Y \cap B$. 
Let $A=\bigcup_{Y_i = X_i} X_i$. Hence $X = A\cup B$ and $Y = A\cup N$.

Consider the proper map $A\sqcup B\longrightarrow X$ which induces an
isomorphism $A\sqcup B - (A\cap B \sqcup  A\cap B) \cong X - (A\cap
B)$. As $\dim A\cap B < n$ we get 
$\Psi\{A\sqcup B, A\cap B \sqcup  A\cap B\}= 
\Psi\{X, A\cap B\}$, hence $\Psi\{A\sqcup B\} = \Psi\{X\} +
\Psi\{A\cap B\}$.

Analogously we get $\Psi\{A\sqcup N\} = \Psi\{Y\} + \Psi\{A\cap N\}$.

Subtracting yields $\Psi\{A\sqcup B, A\sqcup N\} = \Psi\{X,Y\}$ (because $A\cap N = A\cap
B \cap Y = A\cap B$). 

On the other hand $\Psi\{A\sqcup B, A\sqcup N\} = \Psi\{B,N\}$, hence
$\Psi\{B,N\}=\Psi\{X,Y\}$. 

The same reasoning yields $\Psi\{B',N'\}=\Psi\{X',Y'\}$ (with
analogous notations).

Restricting $f$ yields a proper map $(B',N')\longrightarrow (B,N)$
which induces an isomorphism $B'-N' \cong B-N$, and as $\dim N, \dim
N' < n$ we are done.
\end{proof}
This finishes the induction step.

As the Gysin relation holds in $\Knullschlbl(\Vark)$, we get an
induced homomorphism
$\overline{\Psi}:\Knullschl(\Vark)\longrightarrow \Knullschlbl(\Vark)$ 
which can inductively be seen to be an inverse of 
$\Knullschlbl(\Vark) \longrightarrow \Knullschl(\Vark)$.
This establishes Proposition~\ref{presbl}.\qed

\bigskip

We now prove Theorem~\ref{prespaar} by defining an inverse 
$\chi:\Knullschl(\Vark) \longrightarrow
\Knull(\Vark)$ on classes of smooth varieties and
checking compatibility with excision for blow-ups and the Gysin relation.
Let $X$ be a smooth connected variety.
Choose a smooth completion $\overline{X} \supset X$ of $X$ with
$D=\overline{X}-X$ a divisor with simple normal crossings.
\setcounter{step}{0}
\begin{step}
The element $\sum (-\tate)^l[D^{(l)}]$ of $\Knull(\Vark)$ is
independent of the choice of the completion of $X$ 
(we denote it by $\chi(X)$).
\end{step}

\begin{proof}
In view of the weak factorization theorem we only have
to compare the result for a completion $X\subset \overline{X}$ and the
completion $X\subset \Bl{Z}{\overline{X}}=:\overline{X}'$, with $Z$ a smooth and
irreducible closed subvariety of $\overline{X}$ disjoint from $X$
which has normal crossings with $D=\overline{X} - X$.

Let $D_0$ be the irreducible component of $D$ containing $Z$ and let
$\{D_i\}_{i\in I}$ be the remaining irreducible components.
Then the irreducible components of $D'$ are $D'_i=\Bl{Z\cap D_i}D_i$
(where $i\in\{0\}\cup I$) and the exceptional divisor $E$ of the
blow-up.
For $K\subset \{0\}\cup I$ we put $D_K:=\bigcap_{j\in K} D_j$ (where
$D_{\emptyset}$ is understood to be $\overline{X}$), 
$D'_K:=\bigcap_{j\in K} D'_j$, $Z_K:=Z\cap D_K$ and $E_K:=E\cap D'_K$.
Denote the codimension of $Z_K$ in $D_K$ by $d_K$. 
As $Z$ has normal crossings with $D$ we get $D'_K=\Bl{Z_K}{D_K}$ with
exceptional divisor $E_K$, hence by Lemma~\ref{fund2} we have
$$[D'_K]-\tate[E_K]= [D_K]-\tate^{d_K}[Z_K].$$

Denote by $E^{(l)}$ the inverse image of $E$ in ${D'}^{(l)}$. 
For
$l=0,\dots,n$ the preceding equality yields
\begin{align*}
\tate^l[{D'}^{(l)}]&=\tate^l\sum_{|K|=l} [D'_K] +
\tate^l\sum_{|K|=l-1}[E_K]\\
&= \tate^l\sum_{|K|=l} ([D_K]+\tate[E_K] - \tate^{d_K}[Z_K]) +
\tate^l\sum_{|K|=l-1}[E_K]\\
&= \tate^l[D^{(l)}] + \tate^{l+1}[E^{(l)}] - 
\sum_{|K|=l} \tate^{l+d_K}[Z_K] +
\tate^l[E^{(l-1)}], 
\end{align*} 
(for $l=0$ the last term is zero).

For $K\subset I$ we get $Z_{K\cup\{0\}}= Z_K$ and hence $d_{K\cup\{0\}} =
d_K-1$, which yields
$\tate^{|K\cup\{0\}|+d_{K\cup\{0\}}}[Z_{K\cup\{0\}}]=\tate^{|K|+d_K}[Z_K]$
and $\sum_{K\subset I\cup\{0\}}(-1)^{|K|}
\tate^{|K|+d_K}[Z_K]=0$.
Taking the alternating sum we conclude
$$\sum (-\tate)^l [{D'}^{(l)}] = \sum (-\tate)^l [D^{(l)}]$$
as claimed.
\end{proof}
\begin{step}
For a smooth $X$ with connected components $X_i$ we set
$\chi(X):=\sum \chi(X_i)$.
\end{step}

\begin{step}
Let $X$ be a smooth variety, 
let $Z\subset X$ be a smooth
closed subvariety. Let $E$ be the exceptional divisor of
$\Bl{Z}{X}\longrightarrow X$. Then $\chi(\Bl{Z}{X}, E) = \chi(X,Z)$.
\end{step}

\begin{proof}
By decomposing $X$ and $Z$ into connected components we see that it
suffices to prove the claim for connected $X$ and $Z$. 
Denote the dimension of $X$ by $n$.
Choose a smooth completion $\overline{X}$ of $X$ such that 
$D=\overline{X}-X$ is a divisor with simple normal crossings and
$\overline{Z}$ is smooth and has normal crossings with $D$.
So in particular $D\cap \overline{Z}\subset \overline{Z}$ is a simple normal crossings
divisor. Then $\Bl{\overline{Z}}{\overline{X}}$ is a smooth completion of 
$\Bl{Z}{X}$ with 
${D'}:= \Bl{\overline{Z}}{\overline{X}}-\Bl{Z}{X}$ a normal
crossings divisor. The exceptional divisor $\overline{E}$ of this
blow-up is a smooth completion of $E$ with $\overline{E} - E = \overline{E}\cap {D'}$
a simple normal crossings divisor. Denote by $E^{(l)}$ the 
preimage of $\overline{E}$ in ${D'}^{(l)}$ and by $Z^{(l)}$ the preimage
of $\overline{Z}$ in $D^{(l)}$.
Now
\begin{align*}
\chi(\Bl{Z}{X}) & = [{D'} ^{(0)}] - \tate [{D'}^{(1)}] 
 + \dots + (-\tate)^n[{D'}^{(n)}],\\
\chi(E) & = [E^{(0)}] 
- \tate [E^{(1)}] 
 + \dots + (-\tate)^n[E^{(n)}],\\
\chi(X) & = [D ^{(0)}] - \tate [D^{(1)}] 
 + \dots + (-\tate)^n[D^{(n)}],\\
\chi(Z) & = [Z^{(0)}] - \tate [Z^{(1)}] 
 + \dots + (-\tate)^n[Z^{(n)}].
\end{align*}
As ${D'}^{(l)} - E^{(l)} \cong D^{(l)}- Z^{(l)}$ we get
 $[{D'}^{(l)}] - [E^{(l)}] = [D^{(l)}]
- [Z^{(l)}]$ in $\Knull(\Vark)$ and consequently
$$\chi(\Bl{Z}{X},E) = \chi(\Bl{Z}{X}) - \chi(E) = \chi(X) - \chi(Z) 
= \chi(X,Z).$$ 
\end{proof}

\begin{step}
For a smooth connected variety $X$ and $D\subset X$ a smooth divisor we have $\chi(X -D) =
\chi(X) -\tate \chi(D)$.
\end{step}

\begin{proof}
We choose  a smooth completion $\overline{X} \supset X$ 
such that $D_1:= (\overline{X} - X) \cup \overline{D}$ 
is a simple normal crossings
divisor, where by $\overline{D}$ we denote the
closure of $D$ in $\overline{X}$. Let
$D_2:= \overline{X} - X$. Denote by $\Delta^{(l)}$ the inverse image
of $\overline{D}$ in $D_2^{(l)}$. Then $\chi(D) = \sum (-\tate)^l
[\Delta^{(l)}]$ 
and on the other hand
$$\chi (X-D) = \sum (-\tate)^l [D_1^{(l)}] = \sum (-\tate)^l [D_2^{(l)}]
 - \tate \sum (-\tate)^l [\Delta^{(l)}].$$
\end{proof}

\begin{step}
For $X$ smooth and connected $\chi(\proj{1}\times X, \{\infty\}\times X) = \tate\chi(X)$.
\end{step}
\begin{proof}
Let $\overline{X}$ be a smooth completion of $X$ with
$D=\overline{X}-X$ a simple normal crossings divisor.  Then
$\proj{1}\times\overline{X}$ is a smooth completion of $\proj{1}\times
X$ such that the complement of $\proj{1}\times
X$ is $\proj{1}\times D$, a simple normal crossings divisor. This yields
$$\chi(\proj{1}\times X, \{\infty\} \times X) = \chi(\proj{1}\times X)
- \chi(X) = \sum (-\tate)^l([(\proj{1}\times D)^{(l)}]-[D^{(l)}])=\tate\chi(X).$$
\end{proof}

Hence $\chi$ factors through $\Knullschl(\Vark)$. We will now denote the
induced mapping by $\chi$.

\begin{step}
Let $X$ be a smooth connected variety. Then $\chic\chi\langle X\rangle =
\langle X\rangle$ in
$\Knullschl(\Vark)$.
\end{step}

\begin{proof}
Let $\overline{X}\supset X$ be a smooth completion with $D:=
\overline{X} - X$ a simple normal crossings divisor. Let $D_1,\dots, D_n$
be the irreducible components of $D$. We proceed by induction on $n$.

If $X$ is complete then obviously  
$\chic\chi\langle X\rangle = \langle X\rangle$. This
establishes the claim for $n=0$.

Let $X=X_0 \subset X_1 \subset \dots \subset X_n=\overline{X}$
be given by $X_i := X_{i+1} - (D_{i+1}\cap X_{i+1})$.  Then
\begin{align*}
\chic\chi\langle X\rangle &= \chic\chi\langle X_1-(D_1\cap
X_1)\rangle\\ 
&= \chic \chi\langle X_1\rangle - \chic
\chi\langle\proj{1}\times (D_1 \cap X_1), \{\infty\} \times (D_1\cap X_1)\rangle
\end{align*}
by the Gysin relation.  Now $X_1\subset \overline{X}$ is a smooth
completion with $\overline{X} - X_1$ a simple normal crossings divisor with
less then $n$ irreducible components. The same holds for $D_1\cap X_1
\subset D_1$ and for $\proj{1} \times (D_1\cap X_1)\subset
\proj{1}\times D_1$.  Hence by induction we get
\begin{align*} 
\chic \chi\langle X_1\rangle  &- \chic \chi\langle \proj{1}\times (D_1 \cap X_1), \{\infty\} \times (D_1\cap
X_1)\rangle \\ &= \chic \chi \langle X_1\rangle  - \chic \chi \langle \proj{1}\times (D_1 \cap X_1)\rangle 
+ \chic \chi \langle \{\infty\} \times (D_1\cap X_1)\rangle \\
&= \langle X_1\rangle  - \langle \proj{1}\times (D_1 \cap X_1)\rangle  + \langle \{\infty\} \times
(D_1\cap X_1)\rangle
= \langle X\rangle  \text{  again by Gysin.}
\end{align*}
\end{proof}

\begin{step}
The mappings $\chic$ and $\chi$ are mutually inverse.
\end{step}
\begin{proof}
Obviously $\chi\chic[X]=[X]$ for $X$ smooth and
complete, and 
we already know that $\chic\chi\langle X\rangle = \langle X\rangle$ 
for $X$ smooth. 
\end{proof}

\begin{step}
For $X$ a smooth and connected variety the image of $\chi\langle X \rangle$ in 
$\M_k$ is equal to $\D_k(\tate^{-\dim X} [X])=\tate^{\dim X} \D_k([X])$.
\end{step}

\begin{proof}
Let $(\overline{X},D)$ be a smooth simple normal crossings completion. 
Then $$\D_k([X]) = \sum (-1)^l \D_k([D^{(l)}]) = \sum(-1)^l\tate^{l-\dim
X}[D^{(l)}].$$ 
\end{proof}

This finishes the proof of Theorem~\ref{prespaar}. \qed

\section{Relative Grothendieck groups of varieties}
\label{rel}
Let $S$ be a (not necessarily irreducible) variety over $k$.  Let
$\Knull(\VarS)$ be the free abelian group on isomorphism classes
$[X]_S$ of varieties $X$ over $S$ where we impose the relations
$[X]_S=[(X-Y)]_S-[Y]_S$ for $Y\subset X$ a closed subvariety (hence
again $[\emptyset]_S = 0$). We call it the \emph{Grothendieck group of
$S$-varieties}.

This group, too, has simple presentations.
\begin{thm}
\label{presrel}
The Grothendieck group of $S$-varieties has the following alternative
presentations: 
\begin{itemize}
\item[(sm)] As the abelian group generated by the isomorphism classes
of $S$-varieties which are smooth over $k$ subject to the relations
$[X]_S=[X-Y]_S+[Y]_S$, where
$X$ is smooth and $Y\subset X$ is a smooth closed subvariety,
\item[(bl)] as the abelian group generated by the isomorphism classes
of $S$-varieties which are smooth over $k$ and proper over $S$ subject 
to the relations $[\emptyset]_S=0$ and $[\Bl{Y}{X}]_S -[E]_S = [X]_S
-[Y]_S$, 
where $X$ is smooth over $k$ and proper over $S$, $Y\subset X$ is a closed smooth
subvariety, $\Bl{Y}{X}$ is the blow-up of $X$ along $Y$ and $E$ is the
exceptional divisor of this blow-up.
\end{itemize}
Moreover, we get the same group if in case (sm) we restrict to
varieties which are quasi-projective over $S$ or if in case (bl) we
restrict to varieties which are projective over $S$. 
We can also restrict to connected varieties in both
presentations.  
\end{thm}

We proceed in the same way as in the absolute case.
For the moment denote by $\Knullsm(\VarS)$
the free abelian group on isomorphism classes 
$[X]_{S,\sm}$ of varieties over $S$ which are smooth over $k$, modulo
$[X]_{S,\sm}=[X-Y]_{S,\sm}+[Y]_{S,\sm}$ for $X\subset Y$ a smooth 
closed subvariety of $X$.
\setcounter{step}{0}

\begin{step}
The group homomorphism $\Knullsm(\VarS)\longrightarrow \Knull(\VarS)$ is
an isomorphism.
\end{step}

\begin{proof}
We use the same stratification argument as in the absolute case.
\end{proof}
Again we can restrict to smooth connected varieties in the
presentation (sm) by decomposing into connected components.

We can also stratify by smooth varieties which are quasi-projective
over $S$. 

Now we establish the presentation (bl). Denote by
$\Knullbl(\VarS)$ 
the free abelian group on isomorphism classes
of varieties $[X]_{S,\bl}$, where $X$ is a variety which is proper over
$S$ and smooth over $k$, modulo
$[X]_{S,\bl}-[Y]_{S,\bl} =
[\Bl{Y}{X}]_{S,\bl} - [E]_{S,\bl}$ for $Y$ a closed
smooth subvariety of $X$ and $E$  the exceptional divisor of the
blowup of $X$ along $Y$ and $[\emptyset]_{S,\bl} = 0$
(blowing up along $Y$ yields then $[X\sqcup
Y]_{S,\bl}=[X]_{S,\bl} + [Y]_{S,\bl}$).

Decomposing into connected components and noting
that the blow-up along a disjoint union is the successive blow-up along
the connected components one sees that this can also be described as
the free abelian group on isomorphism classes $[X]_{S,\bl}$ of smooth
complete connected varieties with imposed relations
$[\emptyset]_{S,\bl}=0$ and $[\Bl{Y}{X}]_{S,\bl} -[E]_{S,\bl} = [X]_{S,\bl}
-[Y]_{S,\bl}$, where $Y\subset X$ is a closed connected smooth
subvariety. 

\begin{step}
The group homomorphism $\Knullbl(\VarS)\longrightarrow\Knull(\VarS)$
which maps $[X]_{S,\bl}$ to $[X]_S$ is an
isomorphism. 
\end{step}

\begin{proof}
For an arbitrary variety $X$ over $S$ 
there exist a proper
$\overline{X}\rightarrow S$ and an open dense immersion
$X\hookrightarrow \overline{X}$ over $S$. If $X$ is smooth over $k$
and connected we
can even find $\overline{X}$ smooth with $\overline{X}-X$ a simple
normal crossings divisor (by resolution of singularities and
principalization).
As the weak factorization theorem also works over a base variety (note
that a birational map between two irreducible varieties which are
proper over $S$ is automatically proper) 
we can construct an inverse in the same way as in the absolute case.
\end{proof}

Using the fact that we can restrict to generators which are
quasi-projective over $S$ in 
the presentation (sm) and that a connected smooth
variety which is quasi-projective over $S$ has a smooth simple normal
crossings completion which is
projective over $S$ we see that we
can restrict to generators projective over $S$ in the description (bl).

\section{The six operations}
\label{operation}
The product of varieties makes $\Knull(\Var{S})$ a
$\Knull(\Vark)$-module. 
Let $\M_k$ denote
$\Knull(\Vark)[\tate^{-1}]$, let $\M_S$ denote the localization
$\Knull(\VarS)[\tate^{-1}]$. Denote the class $[S]_S$ by $\eins_S$.

\begin{defi}
If $f:S\longrightarrow S'$ is a morphism of $k$-varieties, composition
with $f$ yields a $\Knull(\Vark)$-linear
mapping $f_!:\Knull(\Var{S})\longrightarrow \Knull(\Var{S'})$, hence
we get an $\M_k$-linear mapping $f_!:\M_S\longrightarrow \M_{S'}$.
Pulling back along $f$ yields a $\Knull(\Vark)$-linear mapping
$f^*:\Knull(\Var{S'}) \longrightarrow \Knull(\Var{S})$ and hence an
$\M_k$-linear mapping $\M_{S'}\longrightarrow \M_S$.
\end{defi}

Taking products yields a $\Knull(\Vark)$-bilinear associative exterior product
$$\boxtimes:\Knull(\Var{S})\times\Knull(\Var{T})\longrightarrow
\Knull(\Var{S\times T})$$
and hence an $\M_k$-bilinear associative map 
$$\boxtimes:\M_S\times \M_T\longrightarrow \M_{S\times T}.$$

For morphisms $f:S\longrightarrow S'$ and $g:T\longrightarrow T'$ 
of varieties we get the identities $(f\times g)_!(A\boxtimes B)=f_!(A)\boxtimes
g_!(B)$ and $(f\times g)^*(C\boxtimes D) = f^*(C)\boxtimes g^*(D)$.
If $p:S\times S' \longrightarrow S$ is the projection to the first
factor we get $p^*(A)= A\boxtimes \eins_{S'}$ for $A\in M_S$.
\begin{defi}
Taking fiber products yields an internal $\Knull(\Vark)$-bilinear
symmetric associative product 
$\otimes:\Knull(\VarS)\times\Knull(\VarS)\longrightarrow
\Knull(\VarS)$ and hence an $\M_k$-bilinear associative map
$\otimes:\M_S\times \M_S\longrightarrow \M_S$
\end{defi}
This satisfies $f^*(C\otimes D) = f^*(C)\otimes f^*(D)$ and
$\eins_S\otimes A = A \otimes \eins_S = A$.

For $A\in \M_{S'}$ and $B\in\M_S$ we have $f_!(f^*A\otimes B) =
A\otimes f_! B$. In other words, if we regard $\M_S$ as a 
$\M_{S'}$-module via $f^*$, then $f_!$ is a $\M_{S'}$-module
homomorphism.

For $A,B\in \M_S$ and $C,D\in\M_T$ we get
$$(A\boxtimes C)\otimes(B\boxtimes D) = (A\otimes B)\boxtimes(C\otimes
D),$$ hence the exterior product provides $\M_{S\times T}$ 
with the structure of an
$\M_S\otimes \M_T$-algebra.

\begin{defi} The \emph{duality involution relative to $S$} is the map
defined as follows:

There is a morphism $\D_S:\Knull(\VarS)\longrightarrow \M_S$
which sends a generator $[X]_S$ with $X$ connected and smooth over
$k$, proper over $S$, to  $\tate^{-\dim X}[X]_S$. A relative
version of lemma~\ref{fund2} shows that this is indeed compatible with
blow-up relations.
For $A \in \Knull(\Vark)$ and $B\in \Knull(\VarS)$ we get
$\D_S(AB)=\D_k(A)\D_S(B)$ (this is easily checked on smooth 
proper generators), 
so $\D_S$ can be extended to a $\D_k$-linear morphism
$\D_S:\M_S\longrightarrow \M_S.$
\end{defi}
Indeed $\D_S$ is an involution.
We have already noted in section \ref{univpaar} 
that $\D_k(\tate)=\tate^{-1}$ and that $\D_k$ is a lift of the duality
morphism on the Grothendieck ring of Chow motives.
For $A\in \M_S$ and $B\in \M_T$ we get 
$$\D_{S\times T}(A\boxtimes B) = \D_S(A)\boxtimes\D_T(B).$$ 

\begin{defi}
For $f:S\longrightarrow S'$ define $f^!:=\D_S f^* \D_{S'}$ and 
$f_*=\D_{S'} f_! \D_S$.
\end{defi}
These are both $\M_k$-linear mappings.
\begin{rem}
\label{glatt}
If $f$ is a proper mapping then $f_*$ coincides with $f_!$.
If $f$ is an open embedding then $f^!$ coincides with $f^*$. More
generally, if $f$ is smooth of relative dimension $m$, then for 
$A\in \M_{S'}$ 
we have $f^*A = \tate^m f^! A$.
\end{rem}

\begin{prop}
Given a cartesian diagram as follows
$$\xymatrix{T \ar[d]_\pi\ar[r]^g & T'\ar[d]^{\pi'}\\
              S \ar[r]^f& S'}$$ 
we have
$g_!\pi^* = {\pi'}^*f_!$ and $g_*\pi^! = {\pi'}^!f_*$.
If $f$ is proper or $\pi'$ is smooth, then $g_*\pi^* = {\pi'}^*f_*$.
\end{prop}
\begin{proof}
The only thing to prove is $g_*\pi^* = {\pi'}^*f_*$ for $\pi'$ smooth.
So suppose $\pi'$ is smooth of relative dimension $m$, hence the same
holds for $\pi$.

Let $A\in\M_S$. Then due to Remark~\ref{glatt}
$${\pi'}^*f_*A = \tate^m{\pi'}^!f_* A = \tate^m g_*\pi^! A
g_*(\tate^m\pi^! A) = g_*\pi^* A.$$
\end{proof}
\begin{defi}
Define the \emph{dualizing element} of $\M_S$ by $\DG_S:=\D_S(\eins_S)$.
For $A,B\in \M_S$ define
$$\Hom(A,B) := \D_S(A\otimes \D_S B ).$$
\end{defi}
This yields a mapping 
$\M_S\times \M_S \longrightarrow \M_S$
which is $\D_k$-linear in the first and $\M_k$-linear in the second component.
For $A\in \M_S$ we get $\Hom(A,\DG_S)=\D_S(A)$.

\begin{listprop}
Let $f:S\longrightarrow S'$ be a morphism of $k$-varieties. 
Then the following identities hold:
\begin{align*}
\Hom(A,f_*B)&=f_*\Hom(f^*A,B) \text{ for }A\in \M_{S'} \text{ and }B\in
\M_S,\\
\Hom(f^*A,f^!B)&= f^!\Hom(A,B) \text{ for }A,B\in \M_{S'}.\,\\
f_*\Hom(A,f^!B)&=\Hom(f_!A,B)\text{ for } A\in \M_S\text{ and }B\in
\M_{S'}.
\end{align*}
For $A,B,C\in \M_S$ we get $$\Hom(A, \Hom(B,C)) = \Hom(A\otimes B, C).$$
For $p:S\times S' \longrightarrow S$ the projection to the first
factor and for $A,B\in \M_S$ we have 
$$\Hom(p^*A,p^*B)=p^*\Hom(A,B).$$
\end{listprop}

\begin{proof} 
For we have
\begin{align*}
f_*\Hom(f^*A,B)&= f_*\D_S(f^*A\otimes \D_S B) =
\D_{S'}f_!(f^*A\otimes\D_S B) 
=\D_{S'}(A\otimes f_!\D_S B)\\
&= \D_{S'}(A\otimes\D_S f_* B)=\Hom(A,f_*
B),\\
\Hom(f^*A,f^! B) &= \D_S(f^*A\otimes \D_Sf^!B) = 
\D_S(f^*A\otimes f^*\D_{S'}B)
=\D_S f^*(A\otimes \D_{S'}B)\\ 
&= f^!\D_{S'}(A\otimes \D_{S'}B)=f^!\Hom(A,B),\\
f_*\Hom(A,f^!B)&=f_*\D_S(A\otimes \D_S f^! B)=
f_*\D_S(A\otimes f^*\D_{S'} B)\\
&= \D_{S'}f_!(A\otimes f^*\D_{S'} B) =
\D_{S'}(f_!A\otimes \D_{S'}B) = \Hom(f_!A,B).
\end{align*}
Furthermore,
\begin{align*}
\Hom(A, \Hom(B,C))&= \Hom(A,\D_S(B\otimes\D_S C))\\ &=
\D_S(A\otimes(B\otimes \D_S C))  = \Hom(A\otimes B, C).
\end{align*}
Finally, using $\D_{S\times S'}(C\boxtimes D) = \D_S(C)\boxtimes
\D_{S'}(D)$ for $C\in \M_S$ and $D\in \M_{S'}$ we get
\begin{align*}
\Hom(p^*A,p^*B)
&=\D_{S\times S'}((A\boxtimes\eins_{S'})
\otimes \D_{S\times S'}(B\boxtimes\eins_{S'}))\\
&=\D_{S\times S'}((A\boxtimes\eins_{S'})\otimes (\D_S
B\boxtimes\DG_{S'})\\
&=\D_{S\times S'}((A\otimes\D_S B)\boxtimes \DG_{S'})
= \D_S(A\otimes \D_S B)\boxtimes \D_{S'}(\DG_{S'})\\
&=\Hom(A,B)\boxtimes\eins_{S'} = p^*\Hom(A,B).
\end{align*}
\end{proof}

\section{Equivariant Grothendieck groups of varieties}
\label{equ}
Let $G$ be a finite group. There are also $G$-equivariant versions of 
Grothendieck groups of varieties:

Let $\KnullGstr(\Vark)$ be the free abelian group on $G$-isomorphism 
classes $[G \circlearrowright X]$ (or shorthand $[X]$) of varieties
$X$ 
with a good $G$-action, divided out by
the relations $[X]=[X-Y]+[Y]$ for $Y\subset X$ a closed $G$-invariant 
subvariety. The product of varieties (with diagonal $G$-action) induces
a commutative ring structure on $\KnullGstr(\Vark)$.
Now define $\KnullG(\Vark)$ to be $\KnullGstr(\Vark)$
modulo the ideal $I^G$ generated by expressions of the form
$[G \circlearrowright \Proj(V)]-[G\thickspace\mathrm{trivial} \circlearrowright
\Proj(V)]$, 
where
$V$ is a $k$-vector space on which $G$ acts linearly, $\Proj(V)$
denotes its projectivization and 
$[G\thickspace\mathrm{trivial} \circlearrowright \Proj(V)]$ denotes the class of
$\Proj(V)$ 
with the \emph{trivial}
$G$-action (this is motivated by the fact that on cohomology such an
operation is not visible).
Note that as the class $[G\circlearrowright V]$ is equal to
$[G\circlearrowright \Proj(V\oplus k)]-[G\circlearrowright \Proj(V)]$
this implies that 
$[G \circlearrowright V]-[G\thickspace\mathrm{trivial} \circlearrowright V]$
in $\KnullG(\Vark)$.
\begin{lem}
\label{presequ}
The group $\KnullGstr(\Vark)$ has the following alternative
presentations: 
\begin{itemize}
\item[(sm)] As the abelian group generated by the isomorphism classes
of smooth varieties with good $G$-action 
subject to the relations
$[X]=[X-Y]+[Y]$, where
$X$ is smooth with good $G$-action and $Y\subset X$ is a smooth 
$G$-invariant closed subvariety,
\item[(bl)] as the abelian group generated by the isomorphism classes
of smooth complete varieties with good $G$-action 
subject 
to the relations $[\emptyset]=0$ and $[\Bl{Y}{X}] -[E] = [X]
-[Y]$, 
where $X$ is a smooth complete variety with good $G$-action, 
$Y\subset X$ is a closed smooth $G$-invariant
subvariety, $\Bl{Y}{X}$ is the blow-up of $X$ along $Y$ and $E$ is the
exceptional divisor of this blow-up. 
\end{itemize}
Moreover, we get the same group if in case (sm) we restrict to
quasi-projective varieties or if in case (bl) we
restrict to projective varieties.
We can also restrict to varieties such that $G$ acts transitively on
the connected components in both
presentations.  
\end{lem}
 
\begin{proof}
Stratifying by smooth $G$-invariant equidimensional varieties
establishes the presentation (sm). We can also stratify by
quasi-projective $G$-invariant varieties. We can decompose further
into smooth varieties such that $G$ acts transitively on the connected
components. 

Now note that a smooth connected variety $X$ with a good $G$-action has a good
smooth completion with as complement a $G$-invariant simple normal crossings
divisor:
Take a completion $\overline{X/G}$ of $X/G$ and take the normalization
$\widetilde{X}$ of $\overline{X/G}$ in $X$. This is
finite over $\overline{X/G}$, so the induced $G$-operation is
automatically good. Furthermore
$\widetilde{X}$ is complete and contains $X$ as an open dense subvariety.
An equivariant resolution of singularities of $\widetilde{X}$
making the complement of $X$ a simple normal crossings divisor
then gives a completion
$\overline{X}$ with the desired properties.
Actually as in
remark~\ref{trans} it suffices that $G$ acts transitively on the
connected components of $X$.
The same reasoning as before then establishes the presentation (bl).
Note that in the equivariant version of the factorization theorem the 
centers of the blow-ups need not be connected, but we
can assume that $G$ acts transitively on the connected components.
Restricting to quasi-projective varieties in the presentation (sm)
shows that we can restrict to projective varieties in (bl) (here we
can take any smooth projective equivariant simple normal crossings
completion, as on a projective variety a $G$-action is automatically good).
\end{proof}

Summarizing we get:
\begin{thm}
The ring $\KnullG(\Vark)$ is the free abelian group on smooth
projective (respectively, complete) varieties with good $G$-action
(transitive on the connected components),
with the product described above,
modulo blow-up relations and the ideal $I^G$ generated by expressions 
of the form 
$[G \circlearrowright \Proj(V)]-[G\thickspace\mathrm{trivial} \circlearrowright
\Proj(V)]$.
\end{thm}

We can also consider varieties over a base variety $S$ with good
$G$-action over $S$ (by this we mean that $X\rightarrow S$ is
$G$-equivariant where we give $S$ the trivial $G$-action).
Denote by $\KnullGstr(\VarS)$ 
the free group on isomorphism classes $[X]_S$ of those modulo relations for
closed subvarieties. It has the structure of a
$\KnullGstr(\Vark)$-module provided by the product of varieties with
diagonal $G$-action. 
We define $\KnullG(\VarS)$ to be $\KnullGstr(\VarS)$ modulo the
submodule
$I^G\KnullGstr(\VarS)$ (so in particular it is a $\KnullG(\Vark)$-module).
By the same reasoning as above we get:

\begin{lem}
\label{presrelequ}
The group $\KnullGstr(\VarS)$
has the following alternative
presentations: 
\begin{itemize}
\item[(sm)] As the abelian group generated by the isomorphism classes
of $S$-varieties, smooth over $k$, with a good $G$-action over $S$,
subject to the relations
$[X]_S=[X-Y]_S+[Y]_S$, where
$X$ is smooth with a good $G$-action over $S$ and $Y\subset X$ is a smooth closed $G$-invariant subvariety,
\item[(bl)] as the abelian group generated by the isomorphism classes
of $S$-varieties with good $G$-action over $S$ which are smooth over
$k$  and proper over $S$ subject 
to the relations $[\emptyset]_S=0$ and $[\Bl{Y}{X}]_S -[E]_S = [X]_S
-[Y]_S$, 
where $X$ is smooth over $k$ and proper over $S$ and carries a good
$G$-action over $S$, 
$Y\subset X$ is a closed smooth $G$-invariant 
subvariety, $\Bl{Y}{X}$ is the blow-up of $X$ along $Y$ and $E$ is the
exceptional divisor of this blow-up.
\end{itemize}
Moreover, we get the same group if in case (sm) we restrict to
varieties which are quasi-projective over $S$ or if in case (bl) we
restrict to varieties which are projective over $S$. 
We can also restrict to varieties such that $G$ acts transitively on
the connected components in both
presentations.  
\end{lem}

\begin{thm}
The group $\KnullG(\VarS)$ is the free abelian group on smooth
varieties, projective (respectively, proper) over $S$ with (good) 
$G$-action over $S$ (transitive on the connected components),  
modulo blow-up relations, divided out by the $\KnullGstr(\Vark)$-submodule 
$I^G\KnullGstr(\VarS)$. 
\end{thm}

\end{document}